\documentclass{amsart}

\title{A potentialist conception of ultrafinitism}
\usepackage[pdfauthor={Joel David Hamkins},
    pdftitle={A potentialist conception of ultrafinitism},
    hidelinks]{hyperref}

\author{Joel David Hamkins}
\address[Joel David Hamkins]
{O'Hara Professor of Logic, University of Notre Dame, 100 Malloy Hall, Notre Dame, IN 46556 USA}
\email{jdhamkins@nd.edu}
\urladdr{http://jdh.hamkins.org}
\thanks{This article was adapted from the talk of the same title that I gave at the conference on \emph{Ultrafinitism: Physics, Mathematics, Philosophy} at Columbia University in April 2025. Thanks to Justin Clarke-Doane for helpful comments on an earlier draft.}

\usepackage{latexsym,amsfonts,amsmath,amssymb,mathrsfs}
\usepackage[dvipsnames,svgnames]{xcolor}
\usepackage{wrapfig} 
\usepackage{caption} 
\DeclareCaptionStyle{compact}{font=footnotesize,format=plain,name={Fig},labelsep=period,skip=1.5ex,margin=0pt,oneside,justification=centering}
\DeclareCaptionStyle{leftside}{style=compact,margin={0pt,9pt}}
\DeclareCaptionStyle{rightside}{style=compact,margin={9pt,0pt}}
\usepackage{tikz}
\usetikzlibrary{arrows,arrows.meta,petri,topaths,positioning,shapes,shapes.misc,patterns,calc,decorations.pathreplacing,hobby}
\usepackage[shortlabels]{enumitem} 

\newtheorem{theorem}{Theorem}
\newtheorem*{theorem*}{Theorem}

\newtheorem*{maintheorem*}{Main Theorem}
\newtheorem*{maintheorems*}{Main Theorems}

\newtheorem*{corollary*}{Corollary}
\newtheorem*{corollaries*}{Corollaries}

\theoremstyle{definition}

\newtheorem*{definition*}{Definition}

\newtheorem*{question*}{Question}

\newtheorem*{questions*}{Questions}
\newtheorem*{mainquestion*}{Main Question} 
\newtheorem*{openquestion*}{Open Question} 
\theoremstyle{remark}

\theoremstyle{plain}
\newcommand{\QED}{\end{proof}}

\def\proclaim[#1]{{\bf #1}}
\def\BF#1.{{\bf #1.}}

\def\says#1:#2\par{\item[#1] #2\par}

\newcommand{\Lowenheim}{L\"owenheim}

\newcommand{\C}{{\mathbb C}}

\newcommand{\N}{{\mathbb N}}

\newcommand{\Z}{{\mathbb Z}}
\newcommand{\R}{{\mathbb R}}

\makeatletter
\newcommand{\dotminus}{\mathbin{\text{\@dotminus}}}
\newcommand{\@dotminus}{%
  \ooalign{\hidewidth\raise1ex\hbox{.}\hidewidth\cr$\m@th-$\cr}%
}
\makeatother

\newcommand{\of}{\subseteq}

\newcommand{\set}[1]{\{\,{#1}\,\}}

\newcommand{\Con}{\mathop{{\rm Con}}}

\newcommand{\restrict}{\upharpoonright} 
\newcommand{\satisfies}{\models}

\DeclareMathOperator{\possible}{\text{\tikz[scale=.6ex/1cm,baseline=-.6ex,rotate=45,line width=.1ex]{\draw (-1,-1) rectangle (1,1);}}}
\DeclareMathOperator{\necessary}{\text{\tikz[scale=.6ex/1cm,baseline=-.6ex,line width=.1ex]{\draw (-1,-1) rectangle (1,1);}}}

\newcommand\dbrace{\hskip-1.5em\raise3pt\hbox{\rotatebox[origin=c]{-35}{$\left.\strut^{\phantom{|}}\right\}$}}}

\newcommand\UParroW{{\setbox0\hbox{$\Uparrow$}\rlap{\hbox to \wd0{\hss$\mid$\hss}}\box0}}

\newcommand{\theoryf}[1]{{\rm #1}}

\renewcommand{\setminus}{\raise.3ex\hbox{\rotatebox{-20}{$-$}}} 

\newcommand{\Union}{\bigcup}

\newcommand{\smalllt}{\mathrel{\mathchoice{\raise2pt\hbox{$\scriptstyle<$}}{\raise1pt\hbox{$\scriptstyle<$}}{\raise0pt\hbox{$\scriptscriptstyle<$}}{\scriptscriptstyle<}}}
\newcommand{\smallleq}{\mathrel{\mathchoice{\raise2pt\hbox{$\scriptstyle\leq$}}{\raise1pt\hbox{$\scriptstyle\leq$}}{\raise1pt\hbox{$\scriptscriptstyle\leq$}}{\scriptscriptstyle\leq}}}

   \def\DHLhksqrt#1#2{%
   \setbox0=\hbox{$#1\sqrt{#2\,}$}\dimen0=\ht0
   \advance\dimen0-0.2\ht0
   \setbox2=\hbox{\vrule height\ht0 depth -\dimen0}%
   {\box0\lower0.4pt\box2}}

\def\[#1]{\mathopen{\lbrack\!\lbrack}#1\mathclose{\rbrack\!\rbrack}}
\newbox\gnBoxA
\newbox\gnBoxB
\newdimen\gnCornerHgt
\setbox\gnBoxA=\hbox{\tiny$\ulcorner$}
\global\gnCornerHgt=\ht\gnBoxA
\newdimen\gnArgHgt
\def\gcode #1{%
\setbox\gnBoxA=\hbox{$#1$}%
\setbox\gnBoxB=\hbox{$\bar #1$}%
\gnArgHgt=\ht\gnBoxB%
\ifnum     \gnArgHgt<\gnCornerHgt \gnArgHgt=0pt%
\else \advance \gnArgHgt by -\gnCornerHgt%
\fi \raise\gnArgHgt\hbox{\tiny$\ulcorner$} \box\gnBoxA %
\raise\gnArgHgt\hbox{\tiny$\urcorner$}}
\newcommand{\UnderTilde}[1]{{\setbox1=\hbox{$#1$}\baselineskip=0pt\vtop{\hbox{$#1$}\hbox to\wd1{\hfil$\sim$\hfil}}}{}}
\newcommand{\Undertilde}[1]{{\setbox1=\hbox{$#1$}\baselineskip=0pt\vtop{\hbox{$#1$}\hbox to\wd1{\hfil$\scriptstyle\sim$\hfil}}}{}}
\newcommand{\undertilde}[1]{{\setbox1=\hbox{$#1$}\baselineskip=0pt\vtop{\hbox{$#1$}\hbox to\wd1{\hfil$\scriptscriptstyle\sim$\hfil}}}{}}
\newcommand{\UnderdTilde}[1]{{\setbox1=\hbox{$#1$}\baselineskip=0pt\vtop{\hbox{$#1$}\hbox to\wd1{\hfil$\approx$\hfil}}}{}}
\newcommand{\Underdtilde}[1]{{\setbox1=\hbox{$#1$}\baselineskip=0pt\vtop{\hbox{$#1$}\hbox to\wd1{\hfil\scriptsize$\approx$\hfil}}}{}}

\def\<#1>{\left\langle#1\right\rangle}

\newcommand{\HC}{\mathop{{\rm HC}}}

\newcommand{\ZFC}{{\rm ZFC}}

\newcommand{\ZFCm}{\ZFC^-}

\newcommand{\KM}{{\rm KM}}

\newcommand{\PA}{{\rm PA}}
\newcommand{\PRA}{{\rm PRA}}

\newcommand{\FA}{{\rm FA}}

\newcommand{\cell}[1]{\boxit{\hbox to 17pt{\strut\hfil$#1$\hfil}}}
\newcommand{\head}[2]{\lower2pt\vbox{\hbox{\strut\footnotesize\it\hskip3pt#2}\boxit{\cell#1}}}
\newcommand{\boxit}[1]{\setbox4=\hbox{\kern2pt#1\kern2pt}\hbox{\vrule\vbox{\hrule\kern2pt\box4\kern2pt\hrule}\vrule}}
\newcommand{\Col}[3]{\hbox{\vbox{\baselineskip=0pt\parskip=0pt\cell#1\cell#2\cell#3}}}
\newcommand{\tapenames}{\raise 5pt\vbox to .7in{\hbox to .8in{\it\hfill input: \strut}\vfill\hbox to
.8in{\it\hfill scratch: \strut}\vfill\hbox to .8in{\it\hfill output: \strut}}}
\newcommand{\Head}[4]{\lower2pt\vbox{\hbox to25pt{\strut\footnotesize\it\hfill#4\hfill}\boxit{\Col#1#2#3}}}
\newcommand{\Dots}{\raise 5pt\vbox to .7in{\hbox{\ $\cdots$\strut}\vfill\hbox{\ $\cdots$\strut}\vfill\hbox{\
$\cdots$\strut}}}
\usepackage[backend=bibtex,style=alphabetic,maxbibnames=15,maxcitenames=6,dateabbrev=false]{biblatex}
\renewcommand{\UrlFont}{} 
\renewbibmacro{in:}{\ifentrytype{article}{}{\printtext{\bibstring{in}\intitlepunct}}} 
\DeclareFieldFormat{url}{\UrlFont\url{#1}} 
\DeclareFieldFormat{urldate}{
  (version \thefield{urlday}\addspace%
  \mkbibmonth{\thefield{urlmonth}}\addspace%
  \thefield{urlyear}\isdot)}
\DeclareFieldFormat{eprint:arxiv}{
  \ifhyperref
    {\href{http://arxiv.org/abs/#1}{%
        arXiv\addcolon\nolinkurl{#1}}\iffieldundef{eprintclass}{}{\UrlFont{\mkbibbrackets{\thefield{eprintclass}}}}}
    {arXiv\addcolon\nolinkurl{#1}\iffieldundef{eprintclass}{}{\UrlFont{\mkbibbrackets{\thefield{eprintclass}}}}}}
   \def\DHLhksqrt#1#2{%
   \setbox0=\hbox{$#1\sqrt{#2\,}$}\dimen0=\ht0
   \advance\dimen0-0.2\ht0
   \setbox2=\hbox{\vrule height\ht0 depth -\dimen0}%
   {\box0\lower0.4pt\box2}}

\begin{document}

\begin{abstract}
I shall explore various senses in which ultrafinitism can be fruitfully understood as engaging with a potentialist perspective in mathematics. First, I explain that every model $M$ of the theory of finite arithmetic---arithmetic with a largest number, in which addition and multiplication are merely partial functions---is bi-interpretable with a strictly taller model $M^+$, in which the arithmetic operations on objects taken from the original base model $M$ are totally defined in the extended world $M^+$. More generally, I explain how ultrafinitist ideas emerge in the modal potentialist system consisting of all models of arithmetic under end-extension.
\end{abstract}

\maketitle

\section{Ultrafinitism}

Ultrafinitism is the philosophical view that only comparatively small or accessible numbers exist. According to ultrafinitism, the various extremely large numbers mathematicians conventionally take themselves to describe, such as $10^{100}$ or $10^{10^{10}}$, do not actually exist, and it is a kind of illusion to speak of them. Indeed, often we find ourselves unable to answer basic questions about such numbers. For example, is Skewe's number, the exponential tower $e^{e^{e^{79}}}$\hskip-1em,\hskip1em an integer? We don't know.

Merely denying the existence of all extremely large numbers, however, does not seem to make one an ultrafinitist. After all, other positions in the philosophy of mathematical ontology---I am thinking of certain forms of formalism, fictionalism, nominalism, and so forth---deny in a sense the existence of numbers altogether, but one would not ordinarily classify these views as automatically ultrafinitist. Ultrafinitism, rather, is specifically about a difference in the nature of existence of small versus large numbers. An ultrafinitist accepts the existence of the small or feasible numbers as unproblematic, yet denies the existence of very large numbers. Nevertheless, to be sure, there does seem to be a friendly affinity or overlap between those positions I mentioned and the ultrafinitist attitude toward very large numbers.

Harvey Friedman \cite[p.~4--5]{Friedman2002:Philosophical-problems-in-logic} raised the ``draw the line'' objection with ultrafinitist Alexander Yessenin-Volpin, concerning existence of $2^{100}$.
\begin{quote}\footnotesize
I have seen some ultrafinitists go so far as to challenge the existence of $2^{100}$ as a natural number, in the sense of there being a series of ``points'' of that length. There is the obvious ``draw the line'' objection, asking where in
 $$2^1, 2^2, 2^3,\ldots,2^{100}$$
do we stop having ``Platonistic reality''? Here this \ldots\ is totally innocent, in that it can be easily be replaced by
100 items (names) separated by commas.

I raised just this objection with the (extreme) ultrafinitist Yessenin-Volpin during a lecture of his. He asked me to be more specific.
I then proceeded to start with $2^1$ and asked him whether this is ``real'' or something to that effect. He virtually immediately said yes. Then I asked about $2^2$, and he again said yes, but with a perceptible delay. Then $2^3$, and yes, but with more delay. This continued for a couple of more times, till it was obvious how he was handling this objection. Sure, he was prepared to always answer yes, but he was going to take $2^{100}$ times as long to answer yes to $2^{100}$ then he would to answering $2^1$. There is no way that I could get very far with this.
\end{quote}

The anecdote illustrates the core ultrafinitist idea that number existence becomes increasingly in question as numbers get larger. Larger numbers become inaccessible and subject to feasibility issues. To my way of thinking, the anecdote also opens the door to what I see as an underlying potentialist nature to ultrafinitism, the central theme of this essay. I shall argue, specifically, that ultrafinitism can be fruitfully understood as implementing a potentialist attitude toward number existence---like Yessenin-Volpin in the exchange with Friedman, an ultrafinitist may admit to or accept additional possible numbers into existence if we should simply wait longer. Furthermore, I argue, we shall find ourselves able to express several important ultrafinitist attitudes by adopting a potentialist modal framework, which will also enable us to express and adjudicate disputed aspects of ultrafinitism that might otherwise not even be noticed.

Let me say also that the Yessenin-Volpin/Friedman exchange also happens to illustrate another unfortunate aspect of the ultrafinitism discourse, namely, what I see as a key weakness of many ultrafinitist accounts, which is the general lack of a clear, satisfactory formal theory.

\section{Ultrafinitism in the metatheory}

One of the central philosophical difficulties of ultrafinitism, as I see it, is the challenge it poses for providing a satisfactory formal theory. What is ultrafinitism exactly? In some cases, in practice, ultrafinitist ideas are used not only in the object theory but also in effect employed in the metatheory as well. In some accounts of ultrafinitism, for example, including some accounts given in the talks at the Columbia conference, a formal theory of what are called \emph{feasible} numbers is put forth, with axioms asserting that $0$ is feasible, that the feasible numbers are closed under successor, and also that $2^{100}$ is not feasible.

This is obviously inconsistent, if one could apply the successor operation that many times, and for this reason these accounts of ultrafinitism are often simply a nonstarter amongst mathematicians who encounter them. Meanwhile, there are subtle logical treatments introducing maneuvers to protect against and prevent the inconsistency outcome. For example, one may set up the logical discourse to distinguish the feasible terms, assertions, and proofs from the infeasible instances, in effect allowing semantical judgements only in the feasible case, and in this way one does find a meaningful theory. To my way of thinking, these maneuvers are naturally seen as striving to implement a measure of ultrafinitism in the metatheory. Some accounts of ultrafinitism are quite explicit about having only feasible terms and feasible proofs available in the metatheory. Obviously in an ultrafinitist theory one cannot allow an ordinary approach to meta-theoretic induction on vast terms $1+\cdots+1$, if these are to be taken as meaningful and as large as expected.

The troubling issue with this way of proceeding, of course, is that it takes one outside the usual process of introducing a theory by presenting axioms in a formal language to be considered in classical first-order logic. To require an ultrafinitist metatheory when putting forth a formal theory of ultrafinitism creates a kind of philosophical incoherency in the account. After all, if one takes the metatheory as logically prior to the object theory, then we would seem to need to know already what ultrafinitism is (in the metatheory) before saying what it is (in the object theory). This is a kind of logical circularity, in my view, a confusion that must be sorted in order to provide a satisfactory theoretical account of ultrafinitism. When specifying a formal theory of ultrafinitism, one wouldn't be playing by the usual rules if one were to insist on a measure of ultrafinitism in the metatheory.

Meanwhile, the ultrafinitist perspective on this matter does not see it as circularity, but rather as a kind of ontological honesty---the ultrafinitists want to play by their own rules. In light of this, the circularity criticism cuts both ways. After all, the standard classical logic metatheories include already quite a bit of arithmetic, such as the capacity for induction on terms, formulas, and proofs. In this way, the circularity criticism is aimed right back at the classical theory. Indeed, it is a core ultrafinitist criticism that the standard metatheory of classical arithmetic admits absurd powers of term formation, allowing us freely to form the expression $1+1+\cdots+1$ that would supposedly evaluate to $2^{100}$, even though this is a term that we could never actually form in practice, even if we had been writing down millions of symbols every second since the beginning of time. So the classical arithmetic theory presumes we know what it is (in the metatheory), before saying what it is (in the object theory). Ultimately, for these reasons, the circularity criticism seems to hit both classical and ultrafinitist arithmetic, and perhaps we might fairly judge the dispute on this particular issue therefore to be at a standstill.

\section{A hierarchy of very weak theories}

At my talk for the ultrafinitism conference at Columbia University, I took upon myself the challenge of writing down a list of what I personally consider very weak theories for the foundations of mathematics. At the upper end of what I would consider to be a weak theory I would place Zermelo-Fraenkel set theory $\ZFC$ plus any of the smaller large cardinals. For example, \ZFC\ plus a measurable cardinal, in my view, is a weak theory, much weaker than the far stronger theories routinely considered in set theory, which have shed so much light on the nature and structure of the set theoretic universe. Many of the stronger large cardinal hypotheses have consequences down low at the level of sets of reals in terms of structural results in descriptive set theory. For example, if there are sufficient large cardinals, then every set of reals definable in $\<\R,+,\cdot,0,1,<,\Z>$ is Lebesgue measurable, but this is not provable in \ZFC\ alone, and even a measurable cardinal is not enough (although it is sufficient for the relative consistency result). Weaker large cardinal theories would posit merely Ramsey cardinals, or perhaps strongly unfoldable cardinals, or totally indescribable cardinals. Eventually one will get down to a Mahlo cardinal or perhaps a proper class of inaccessible cardinals, a hypothesis that is equivalent to the Zermelo-Gorthendieck axiom of universes commonly used in category theory. Weaker still would be \ZFC\ plus just one inaccessible cardinal and below that is Kelley-Morse \KM\ set theory, a very weak set theory. Weaker still is \ZFC\ plus a proper class of otherworldly cardinals.

Eventually, we shall weaken our theory down to plain \ZFC\ itself, a very weak theory---the pervasive independence phenomenon reveals so much about what \ZFC\ doesn't settle. Nevertheless, one can weaken still further, down to $\ZFCm$, which is set theory without the power set axiom, a theory that holds in the realm of hereditarily countable sets, which satisfies the countabilism axiom, asserting that every set is countable. These sets can therefore be captured in effect by second-order number theory, and this brings us finally into the realm of arithmetic. There is a richly studied hierarchy of second-order number theories in reverse mathematics, all of them extremely weak theories from the perspective of large cardinal set theory. In all the usual theories of reverse mathematics, for example, one can define and express sets in the formalism of the theories, and yet the theory does not prove that these sets exist, which I take as a fundamental weakness of these theories.

Finally, eventually, we give up the second-order arithmetic context and enter the realm of first-order arithmetic with Peano arithmetic \PA---the theory of a discretely ordered semi-ring with the first-order induction scheme---true in the standard model of arithmetic $\<\N,+,\cdot,0,1,<>$. To my way of thinking, this is the true natural realm of finite mathematics and finitary reasoning, although one will find this latter point to be vigorously contested by advocates of weaker foundations. One can restrict the induction scheme, for example, to $\Sigma_n$ formulas, making the theory $I\Sigma_n$---one might take it all the way down to $I\Sigma_1$, which allows induction only for purely existential assertions, sufficient to prove the totality of exponentiation and other recursive definitions. Or one might weaken it to the theory of primitive recursive arithmetic \PRA, which includes the recursive definitions of all primitive recursive functions, but only quantifier-free induction in this language.

The theory of bounded induction I$\Delta_0$ in the theory in the language of arithmetic has induction only for $\Delta_0$ assertions, that is, only for assertions with all quantifiers bounded. This theory proves all the most elementary facts about addition and multiplication, but meanwhile, it does not prove the totality of exponentiation. It is consistent with I$\Delta_0$ that there are numbers $n$ for which $2^n$ does not exist. In this sense, the theory begins to engage with several core ultrafinitist ideas, although not the most radical versions of ultrafinitism, since the existence of $2^{100}$ specifically is provable in I$\Delta_0$---the lack of totality can arrive only much higher, with nonstandard instances. Nevertheless, one might take it as a core ultrafinitist commitment that addition and multiplication are innocent, while exponentiation is not total, and this is a perspective that I$\Delta_0$ implements quite satisfactorily.

One can weaken the induction axiom further, to the quantifier-free formulas, resulting in the theory of open induction, and below this one would find Robinson's theory $Q$, an absurdly weak theory of arithmetic, of scarcely any interest as a foundation of mathematics, in which one cannot prove even that addition is commutative or that multiplication distributes over addition.\footnote{The theory $Q$ is interesting in that it provides a minimal theory that is sufficient to support the conclusions of the incompleteness theorem; but I am flabbergasted by the occasional proposals one hears to consider it somehow as a foundational theory, a role for which it strikes me as obviously unsuited.} Nevertheless, even this ridiculous theory $Q$ proves that $2^{100}$ exists.

\section{The theory of finite arithmetic \FA}

All the theories I have mentioned so far prove that the natural numbers are closed under successor, and the various ultrafinitist theories to which I have alluded introduce concepts of feasibility that aim in effect to block this process from reaching the realm of infeasibility.

Meanwhile, there is a completely different approach to ultrafinitism that gives rise to a totally different genre of theory, which I should like now to discuss in detail. Namely, in this alternative account of ultrafinitism, one denies the closure under successor and instead simply posits outright the existence of a largest number. The idea that there is a largest number was put forth by Doron Zeilberger, for example, in his plenary talk at the Columbia University conference.

Let me introduce and formalize here a specific such theory, what I call the \emph{theory of finite arithmetic}, abbreviated \FA, although this theory is also known as $\PA^{\!\textup{top}}$ in the research literature. This is the theory of arithmetic with a largest number, expressed in the usual language of arithmetic $\set{{+},{\cdot},0,1,N,{<}}$, but augmented with a constant $N$ for the largest number. The addition and multiplication operations are taken as merely partial functions, and remain undefined in effect when the result would be too large. One can formalize partial functions in first-order logic officially by treating them via their graph relations (see \cite{Hamkins2024:Definite-descriptions} for several ways to do this, which all work fine in this context); for simplicity I will use functional notation in this account, but let us think of them as partial functions.

The axioms of finite arithmetic \FA\ are:
\begin{enumerate}
  \item (Order) The $<$ relation is discrete linear order, with least element $0$ and largest number $N$.
  \item (Successor) The successor of any number $n$ is $n+1$, which is defined for every $n<N$.
  \item (Arithmetic) Addition and multiplication are partial functions obeying the following recursive definitions, taken in the sense that if the right-hand side is defined, then so is the left-hand side and they are equal:
 $$n+0=n\qquad\qquad n+(m+1)=(n+m)+1$$
 $$n\cdot0=0\quad\qquad\qquad n\cdot(m+1)=n\cdot m+n$$
  \item (Induction) The full induction scheme, for any formula $\varphi$ in the language of finite arithmetic
$$\left[\varphi(0)\text{ and }\forall n<N\, \bigl(\varphi(n)\to\varphi(n+1)\bigr)\right]\to\forall n\ \varphi(n).\quad $$
\end{enumerate}
In this theory, one can prove the expected basic arithmetic facts, as much as they would be expected to hold below $N$, such as associativity, commutativity, distributivity, and so forth, with proofs by induction largely following the usual development in Peano arithmetic \PA. Since exponentiation will not be available in \FA, however, the opportunities for arithmetization and coding will be far more limited in \FA\ in comparison with \PA.

The canonical models of finite arithmetic \FA, of course, are those arising by truncation of the standard model of arithmetic $\N$ at a number $n\geq 1$, resulting in the model I denote $\N\restrict n$, built on domain $\set{0, 1, 2, \ldots, n}$, with the addition and multiplication operations restricted to instances where the result is at most $n$. It is easy to see that \FA\ is true in any such truncated model $\N\restrict n$, and so \FA\ is part of the common theory of these models.

The theory of finite arithmetic is surely inspired by these truncation models $\N\restrict n$ and naturally they are the principal intended models of this theory. A mathematician might naturally expect that these are the only models of the theory, but we know this cannot be right because the upward \Lowenheim-Skolem theorem will give us nonstandard infinite models of the theory. Nevertheless, a contemplative mathematician might still reasonably expect that the theory of finite arithmetic \FA\ should be exactly the common theory of these models $\N\restrict n$. Perhaps such a mathematician hopes to reason in finite arithmetic by imagining simply that the universe is some standard truncation model $\N\restrict n$, except only that the particular value of $n$ is not known. Regrettably, this view of the situation is not correct and \FA\ is simply not the same as the common theory of the truncations $\N\restrict n$. Indeed, what I claim is that the common theory of all $\N\restrict n$ is actually too complex for this, for it can have no computable axiomatization. As a computational matter, we aren't actually able to reason fully in that common theory. To see this, observe simply that a Turing machine fails to halt if and only if it fails to halt in every truncation model $\N\restrict n$, and so if we had a computable axiomatization of this common theory we could search for a proof and thereby recognize the non-halting programs. Since we can also recognize the halting programs, this would give us a computable solution of the halting problem, which is impossible. Therefore, there can be no computable axiomatization of the common theory of all the truncation models $\N\restrict n$. Similar reasoning shows that the common theory of all $\N\restrict n$ is in fact a $\Pi^0_1$-complete set.

So \FA\ is strictly weaker than the common theory of the standard truncation models $\N\restrict n$, and consequently there must be other strange nonstandard models of finite arithmetic \FA. As mentioned, we get nonstandard models of \FA\ by the \Lowenheim-Skolem theorem, or one can use the usual kind of compactness arguments or the ultraproduct construction---just take the ultraproduct $\prod_n (\N\restrict n)/\mu$ modulo a nonprincipal ultrafilter $\mu$ on $\N$ (although of course this will be a model of the common theory). More generally, every model $M\satisfies\PA$ of Peano arithmetic has all its truncations $M\restrict n$, including nonstandard truncations, being models of \FA. And these can satisfy assertions not true in any truncation $\N\restrict n$, such as $\neg\Con(\PA)$. It is sufficient merely that $M\satisfies\text{I}\Delta_0$ in this argument. That is, what I claim is that for any such model $M$ of bounded induction, all the truncation models $M\restrict n$ will satisfy the full induction scheme, since this corresponds to an instance of bounded induction in $M$, as we can bound all the quantifiers by $n$ in the ambient structure $M$. So these will be models of finite arithmetic \FA.

I find it interesting to observe that the existence of such nonstandard models of finite arithmetic, if one takes this as an ultrafinitist theory, opens up the possibility of ultrafinitist nonstandard analysis, a combination that might otherwise seem ridiculous. But this way of proceeding with it strikes me as  perfectly sensible and fruitful.

\section{The ad hoc criticism}

There is something a little unsettling about taking the formal theory of finite arithmetic \FA\ seriously as one's best account of the fundamental nature of arithmetic. How could this theory be the correct, final philosophical account of the nature of finite numbers? The sharp cut off at the largest number $N$ seems arbitrary---it is \emph{ad hoc}. Why should the numbers stop sharply right at that point and not go on a little further? There would be something deeply contingent for the numbers to stop at that particular number $N$. Why not define $N+1$ somehow and provide a suitable, meaningful extension of the arithmetic operations up to this new number? For this reason, it strikes many as absurd that our best, ultimate theory of arithmetic would ever posit the existence of a largest number.

\section{Interpreting a taller model of \FA\ inside any model of \FA}

I should like to give mathematical legs to the \emph{ad hoc} criticism by explaining exactly how every model of \FA\ can interpret a strictly taller model of this theory. The fact of the matter is that inside any model of \FA\ we can interpret a taller model of \FA\ in which the formerly largest number $N$ now achieves its square $N^2$ or indeed $N^3$ or much more. Indeed, it follows from this that in fact, all sums and products of members of $M$ become fully defined in this taller model, even if they were not meaningful in the original model $M$. That is, all the undefined cases of sums and products in $M$ become meaningful in the taller model. Ultimately, by iterating these interpretations, we shall show that every model of \FA\ arises via truncation from a model of bounded induction I$\Delta_0$.

To begin, let us suppose that we have a model $M\satisfies\FA$ with largest number $N$. It follows by the induction axiom that there is some largest $b$ such that $b^2$ exists in $M$. This number $b$ is something like $\sqrt{N}$, possibly rounded down, and indeed inside the taller model we aim to construct, this will be literally correct.

The main thing to observe about the number $b$ is that addition and multiplication is totally defined in $M$ for all numbers below $b$. We may add and multiply such small numbers freely, although the result may take us out of that smallness realm.


Let me assume for nontriviality that the largest number $N$ inside $M$ is already quite substantial, say, at least one hundred---and we have found the largest number $b$ for which $b^2$ exists. I propose that we work in base $b$ with, say, five-digit numbers \texttt{ABCDE}. All individual digits are less than the base $b$. For convenience, we shall allow leading zeros in the representation, such as with \text{00129}, so that all the number strings will have exactly five digits. We may consider these five-tuples as lists of five numbers from $M$, which we shall view from the perspective of $M$ as representing or interpreting some kind of comparatively huge imaginary number. In this way, the model $M$ will have access in a way to a richer realm of numbers, which may be seen as strictly imaginary or merely possible from the point of $M$, but which nevertheless will constitute a strictly taller realm of numbers, a taller but perfectly coherent model of finite arithmetic, yet fully understandable to someone whose numbers are limited to $M$.

To clarify, we are not performing arithmetization somehow to code these sequences inside $M$, and indeed that would be impossible since $M$ does not have enough numbers, and we are certainly not giving these number strings a numerical value inside $M$, since those values would become too large---the number string \texttt{10000}, for example, is trying to represent $b^4$, which is larger than $N$, the largest available number in $M$. Rather, we are simply considering lists of five numbers from $M$, all less than $b$, to be used as parameters in a formula that might express something about the number string that they represent.

For example, working entirely inside $M$ we can judge the order relation on these number strings, since we place \texttt{ABCDE} before \texttt{PQRST} just in case it appears earlier in the lexical order of the values of those digits, which can be defined by comparing the digits in turn. That is, there is an expressible formula $\text{LessThan}(A,B,C,D,E,P,Q,R,S,T)$ which holds in $M$ when these variables represent suitable number strings in the lexical order.

Next, again whilst working inside $M$, we can define a successor operation on these number strings, that is, the string \texttt{FGHIJ} that would correspond to \texttt{ABCDE}$+$\texttt{1}.
\begingroup
\setlength{\jot}{0pt} 
\begin{align*}
        &\texttt{ABCDE} \\
 +\quad &\texttt{\ \ \ \ 1} \\[-7pt]
   \cline{1-2}
        &\texttt{FGHIJ}
\end{align*}
\endgroup
To do this, we simply perform the usual grade-school addition algorithm to add $1$ to a number, adding $1$ to the digit \texttt{E}, and in the case that \texttt{E} was already the largest digit below $b$, then we take \texttt{J} as \texttt{0} and perform the carry operation, which means adding $1$ to digit \texttt{D} and so forth. This entire process can be undertaken inside $M$ using only the arithmetic resources that are available there. Adding $1$ to the very largest five-digit number \texttt{ZZZZZ}, where $\texttt{Z}=b-1$ is the largest base $b$ digit, of course, would result in a 6-digit number \texttt{100000}, and so that instance of the operation will remain undefined amongst the five-digit numbers.

We may similarly adopt the standard grade-school algorithm for addition, which will allow us to add using these number strings.
\begingroup
\setlength{\jot}{0pt} 
\begin{align*}
        &\texttt{ABCDE} \\
 +\quad &\texttt{FGHIJ} \\[-7pt]
   \cline{1-2}
        &\texttt{PQRST}
\end{align*}
\endgroup
The point is that in order to carry out ordinary addition in base $b$ you only ever need to add numbers less than $b$ at any stage in the process. As elementary school teachers everywhere will attest, if you know your addition tables up to 10, then you can add any two numbers at all in decimal notation using the grade-school addition algorithm. Similarly, whilst staying inside $M$ and using only the arithmetic resources available there, we can add five-digit numbers in base $b$, and also recognize whether the result remains a five-digit number. What I mean formally is that the grade-school algorithm is definable in the language of finite arithmetic, in the sense that there is a definable relation $\text{Plus}(A,B,C,D,E,F,G,H,I,J,P,Q,R,S,T)$ that holds in $M$ just in case the sum would come out like that for those digits correctly. To express this in $M$, we do not need to do any arithmetization coding or indeed any arithmetic at all beyond adding the individual digits, which are all less than $b$, and so this can be done in $M$.\goodbreak

It is just the same with multiplication. We can implement the grade-school multiplication algorithm, which at each step involves adding or multiplying the digits.
\begingroup
\setlength{\jot}{0pt} 
\begin{align*}
        &\texttt{\ \ ABC} \\
 \times\quad &\texttt{\ \ \ DE} \\[-7pt]
   \cline{1-2}
        &\texttt{\ FGHI}\\
  +\quad &\texttt{JKLM}\\[-7pt]
   \cline{1-2}
      &\texttt{NOPQR}\\
\end{align*}
\endgroup
Note that we only care to give an answer when the result can be expressed in our five digit number strings. All this is expressible inside $M$ as a property of the input digits.


This number-string idea ultimately serves, I claim, to interpret a taller model of finite arithmetic inside $M$. To see this, let $M^+$ denote the model having as domain the five-digit base $b$ number strings \texttt{ABCDE} as formed in $M$, with the order $<$ as we defined it and with the addition and multiplication as defined there by the grade-school algorithms. We interpret the largest-number constant symbol $N$ by the string \texttt{ZZZZZ}, where $\texttt{Z}=b-1$, the largest digit available in base $b$. All this structure is expressible and definable inside $M$.

What I claim is that $M^+$ is a model of finite arithmetic \FA. It is clear that the order $<$ we defined on the number strings \texttt{ABCDE} is discrete, and that adding one as it is defined on the number strings is the successor with respect to this order. Furthermore, the arithmetic axioms of addition and multiplication are valid for the same kind of reasons that show the legitimacy of the grade-school algorithms for arithmetic. Finally, I claim that the induction scheme holds in $M^+$. To see this, assume in $M^+$ that $\varphi(0)$ holds and that $\forall x\,\bigl(\varphi(x)\to\varphi(x+1)\bigr)$ for some formula $\varphi$ in the language of finite arithmetic. Inside the ground model $M$, we can express whether $\varphi(x)$ holds or not in $M^+$ of any given five-digit number string \texttt{ABCDE}. And we can look for the smallest first digit \texttt{A}, which may be part of a counterinstance to $\varphi$, and then the smallest second digit \texttt{B}, for which there is a counterinstance starting \texttt{AB}, and so forth. In five steps, applying induction this way only in $M$, we will thereby find a lexically least five-digit counterinstance \texttt{ABCDE} to $\varphi$, which will violate our assumptions that $\varphi$ was inductive in $M^+$. So $M^+$ is a model of \FA.

Next we may observe that $M$ embeds naturally into $M^+$ as an initial segment, since every element of $M$ up to $b^2$ can be represented in the form $bn+k$, where $n,k<b$, which will be represented by the two digit numbers \texttt{000DE}, and then since $(b+1)^2$ does not exist, there are at most $2b$ many additional numbers on top, which can be represented in the form $b^2+bn+k$. And so in $M$ we can define the representation of every number in the form \texttt{00CDE}, where \texttt{C} is either $0$ or $1$ for these two situations, respectively, and thereby embed $M$ as an initial segment of $M^+$.

The model $M^+$ is much taller than $M$, and indeed $N^2$ exists in $M^+$, since \texttt{ZZZZZ} represents $b^5-1$ in $M^+$, and this is larger than $N^2$. The general conclusion is that every model of finite arithmetic, with largest number $N$, interprets a taller model of finite arithmetic in which $N^2$ exists. In particular, it follows that for any numbers $x,y\in M$ in the original model, the sum $x+y$ and product $xy$ are defined in $M^+$.

I should like particularly to emphasize this point. What it shows is that the lack of a meaning for $x+y$ or $xy$ in $M$ was contingent. That is, by simply moving to $M^+$ these arithmetical values, which were undefined and meaningless in $M$, become fully defined and numerically meaningful in $M^+$. Indeed, to foreshadow the potentialist outlook that I shall discuss later, we can say that although the sum and product of numbers in $M$ in some cases may not be actually defined, nevertheless they are possibly defined, in that there is an interpretable larger model of the theory in which that particular sum or product is defined.

A more sophisticated, deeper observation is the fact that the two models $M$ and $M^+$ are bi-interpretable. We have given an interpretation of $M^+$ in $M$ by means of the number strings \texttt{ABCDE}, and furthermore $M$ can define the copy of itself inside $M^+$, which is definable in $M$. Conversely, $M^+$ can interpret $M$ as the relevant truncation and it can see how its numbers all have five-digit representations in base $b$. So $M^+$ can see how it is copied into the interpretation of it arising in the copy of $M$. Each model can see how it arises in the interpreted copy of the other model, and this is what it means to be a bi-interpretation.


My main philosophical point here is that I take this interpretation construction to give mathematical legs to the ad hoc criticism I mentioned earlier. Finite arithmetic \FA\ purports to give us an account of arithmetic by providing a realm of numbers from 0 up to $N$. But by using only those ontological resources, we can form all the five-digit number strings and thereby interpret a much taller system. So there was no need after all for us to have stopped at $N$. There would have been nothing incoherent about having more numbers strictly beyond $N$. For the original model to have stopped at $N$ was an unnecessary, arbitrary move. I take this to show that the concept of arithmetic provided by the original model cannot be the correct, final account of the finite numbers.


Let me observe next the amazing things that happen when we should simply iterate the construction. Begin with any model of finite arithmetic that you like, $M\satisfies\FA$. Inside this model, we can interpret a strictly taller model $M^+$, in which all the sums and products of the former numbers become meaningful. And inside this taller model $M^+$, we can interpret a still taller model $M^{++}$. Since interpretations are transitive, in truth we could have done it in one step, interpreting $M^{++}$ in $M$ simply by using longer number strings. But we may continue iterating indefinitely, interpreting taller and taller models of finite arithmetic, in such a way that the sums and products of numbers in any of the models $M^{+n}$ become defined in the next model $M^{{+}(n+1)}$.
 $$M\of M^+\of M^{++}\of M^{+++}\of\cdots$$
Each model can see itself copied as an initial segment of the next, and by suitable replacing of these copies with the real thing, if this would be desired, we can assume without any loss that each model is an actual substructure of the next. Let us now take the limit structure to which these models are converging.
 $$M^*=\Union_n M^{+n}$$
As one ascends in the tower, moving from $M$ to $M^+$ and to $M^{++}$ and so on, the arithmetic operations of each model are stabilized in the next. Indeed, since the arithmetic operations become successively defined---the sum and product of any two individuals becomes defined in the next model of the tower after those individuals come into existence---the limit model $M^*$ therefore realizes its addition $+$ and multiplication $\cdot$ operations as totally defined. So the limit model $M^*$ is a model of some kind of arithmetic. What kind of arithmetic theory can we expect to hold in this interpretation limit model?

If the original model $M$ was a truncation of the standard model $\N\restrict n_0$, then it will follow that all the $M^{+n}$ are standard and the limit model $M^*$ in this case will simply be a copy of the standard model of arithmetic $\N$, which of course satisfies an excellent arithmetic theory.

In the general case, however, the initial model $M$ may be a nonstandard model of \FA, in which case all the models $M^{+n}$ and $M^*$ will also be nonstandard. And yet, every model $M^{+n}$ along the way was a model of \FA, which includes the full induction scheme in the language of arithmetic. So perhaps one might hope to realize Peano arithmetic \PA\ in the limit model $M^*$, with the full induction scheme. After all, this is what we saw in the standard case, where $M^*$ was the standard model $\N$, which satisfies \PA.

Unfortunately, however, except in the standard case, we will not generally get a model of \PA\ in the limit. One can see this by observing that if $N$ was the largest number in $M$, then we achieve $N^2$ in $M^+$, and $N^4$ in $M^{++}$ and so forth, but in this construction we will only ever achieve $N^k$ for standard numbers $k$. So if $N$ was nonstandard, then there will be no such number as $N^N$ in $M^*$, and for this reason in this case $M^*$ will not be a model of \PA.

In fact, the expectation that one might achieve full induction in $M^*$ was misplaced from the start, the reason being that the full induction scheme in \FA\ that $M^*$ is inheriting corresponds in effect to \emph{bounded} induction, as the unbounded quantifiers in any model $M^{+n}$ are bounded by the largest number of that model. In this way, however, we can see that the interpretation limit model $M^*$ will be a model of bounded induction I$\Delta_0$. If $\varphi$ is a formula in the language of arithmetic with all quantifiers bounded, and we assume $M^*\satisfies \varphi(0)$ and $\forall x\,\bigl(\varphi(x)\to\varphi(x+1)\bigr)$, then these hypothesis will also be true in every model $M^{+n}$ along the way, meaning that every $M^{+n}$ satisfies $\forall x\,\varphi(x)$. And since $\Delta_0$ truth is absolute between $M^*$ and its initial segments $M^{+n}$, this means we conclude $\forall x\,\varphi(x)$ in $M^*$ as well. So the interpretation limit model $M^*$ is a model of bounded induction I$\Delta_0$.

\section{Characterization of the models of finite arithmetic}

The summative general conclusion here is that every model of finite arithmetic \FA\ arises by truncation from a model of bounded induction.

\begin{theorem}[Paris]
 Every model of finite arithmetic $M\satisfies\FA$ arises by truncating a model of bounded induction $M^*\satisfies\text{I}\Delta_0$ at some number $n\in M^*$, possibly nonstandard, and conversely.
\end{theorem}

\noindent That is, the models of \FA\ are exactly the models of the form $M^*\restrict n$ for some $M^*\satisfies\text{I}\Delta_0$. I had proved this result on my own a few years ago, using the argument I have just described, but shortly afterward I came to realize that actually it is a well-known result in the field of models of arithmetic, due originally to Jeff Paris, and one can find accounts of it in \cite[thm.IV.2.2]{HajekPudlak1998:Metamathematics-of-first-order-arithmetic} and \cite[thm.13.11]{Wong2015:Arithmetic-with-a-top}, following a similar proof idea to the base $b$ number string interpretation I provided above.

Let me point out a few further consequences we can extract from the proof. Namely, first, the extension model $M^*$ we constructed is the unique smallest such extension of $M$ to a model of I$\Delta_0$, and it will appear as an initial segment of all other such extension models. This is because every model of I$\Delta_0$ extending $M$ will see the five-digit number strings in $M$ and agree with $M^+$ about how they add and multiply, and then also with the next model $M^{++}$ and so forth. The basic situation is that any other extension $M^{**}$ can run the interpretations of $M^+$, $M^{++}$, and so forth, and the structure and nature of these models is determined completely by $M$, so it will agree with $M^*$ up to that height, although $M^{**}$ might be taller still.

In particular, let me highlight a particular consequence of this, a consequence that is relevant for my discussion of potentialism in the next section. Namely, although the arithmetic operations are not totally defined in the original model $M$, in that some sums and products were not defined there, nevertheless the way that these sums and products will come to be defined in taller models of \FA, and furthermore the basic properties that these sums and products will come to have in such extensions, regarding whether they are prime or square free and so on, is completely determinate in all such extensions. All extensions of $M$ to a taller model $M^+$ of finite arithmetic on which the sums and products of members of $M$ exist will agree on all the $\Delta_0$ properties of those sums and products, and so there is no contingency at all in the basic arithmetic properties that the numbers of $M$ might come to have in a taller model of \FA\ or even I$\Delta_0$.

Another consequence of the proof is the observation that I$\Delta_0$ does not prove the totality of exponentiation. We have observed that if $N$ is the largest number of the original model $M$, and this is nonstandard, then in $M^*$ we shall have only the standard powers $N^k$ for standard numbers $k$. There will be no $N^N$ and therefore no $2^{N^2}$ in the limit model $M^*$. So I$\Delta_0$ does not prove that $2^n$ exists for every number $n$.

Let me also make a philosophical point here. Namely, by realizing every model of finite arithmetic as the truncation of a model of I$\Delta_0$, we come to see that these two visions of ultrafinitism, despite initial appearances, are intimately connected. I had presented \FA\ as a formulation of ultrafinitist ideas---in positing a largest number it denies the existence of extremely large numbers. And I had also discussed the theory of bounded induction I$\Delta_0$ as having an ultrafinitist nature, in that it is closed under successor, yet is unable to prove the existence of exponentials $a^b$. Although I$\Delta_0$ may seem at first to partake far more of an infinitary nature than \FA, nevertheless we have now proved that every model of \FA\ leads to a unique smallest model of I$\Delta_0$ extending it. So every model of finite arithmetic has already sufficient ontological resources to interpret taller such models and then still taller such models again, each of them bi-interpretable with the original model. In light of these bi-interpretations, these taller models are already present in the smaller model, being represented by suitable number strings \texttt{ABCDE}\dots. By taking these successive interpretations altogether, each of which is present in the smaller model in this way, we construct the unique smallest limit model $M^*$ of the theory of bounded induction I$\Delta_0$. In this sense, we would encourage any ultrafinitist who finds the theory of finite arithmetic \FA\ appealing, to consider to extent to which their commitments to \FA\ might also justify I$\Delta_0$, in light of the fact that one can build models of this latter theory on the ontology of the former.

%
%
%
%
%

\section{Arithmetic potentialism}

Finally we come into the main theme of my essay. I should like to discuss what I see as an underlying potentialist nature to ultrafinitism. I propose that we may fruitfully view the ultrafinitist perspective in a potentialist light, which will help illuminate its philosophical commitments, whilst also enabling a formal treatment of various ultrafinitist theories. What is more, I believe that the potentialist perspective brings to light certain fundamental issues on the nature of mathematical existence on which ultrafinitists might disagree, but which are most naturally discussed and adjudicated in a potentialist setting.

At the Columbia University conference, Sam Buss mentioned that Ed Nelson had expressed ideas having a certain affinity with a potentialist outlook, in particular, the idea that things become true in arithmetic as you develop the theory---perhaps the twin primes conjecture could become true or the negation, depending on how the theory develops.

From the point of view of this essay, however, my entry into potentialism arises instead from a semantical perspective regarding the models of arithmetic, including models of the weak or the ultrafinitist theories. We saw, for example, how every model $M\satisfies\FA$ of finite arithmetic extends to taller models $M^+$ and $M^{++}$ and so forth, with which it is bi-interpretable, and $M$ extends ultimately to the limit model $M^*\satisfies\text{I}\Delta_0$ of bounded induction. My view on this is to take it directly as a form of potentialism. Even if a finite-arithmetic ultrafinitist does not agree with the $M^*$ limit construction, nevertheless it seems to be in accordance with ultrafinitism to allow the move from $M$ to $M^+$, and this kind of move already leads exactly to the potentialist picture for arithmetic that I would like to paint.


As I see it, ultrafinitism is centered upon the concept of a realm of arithmetic feasibility. And arithmetic potentialism is all about having realms of feasibility, about how there might in fact be different realms of feasibility at different scales of feasibility. An extremely limited concept of what is feasible, for example, might be expanded, allowing greater access to the resources of time, space, or even computational power, to a broader conception of what is feasible. And that conception might expand still further to larger realms of feasibility.


In this way, we achieve a semantical picture of arithmetic potentialism, one where there is a space of possible worlds or arithmetic realms of possibility, instantiating the idea that you can have more and more, even if at no point do you ever have it all. These various realms of arithmetic can be taken as constituting a Kripke model of possible worlds. Each such world can be seen as a particular instantiation of ultrafinitist ideas for a particular concept of feasibility.

Current philosophical work has thus emphasized the \emph{modal} nature of potentialism. The various universe fragments are possible worlds in a potentialist system, which enables a modal vocabulary and the modal semantics. Namely, an assertion $\varphi$ is \emph{possible} at a world, written $\possible\varphi$, if it holds in some larger world. And $\varphi$ is \emph{necessary} at the world, written $\necessary\varphi$, if it holds in all the larger worlds. The modal vocabulary allows us precisely to express various potentialist facts about the current realm of numbers in which we may find ourselves.


\begin{wrapfigure}{r}{.25\textwidth}\vskip-3ex\hfill
\begin{tikzpicture}[dot/.style={circle,fill,inner sep=1pt}]
\draw[Red,fill=Gold!20] (0,1.6) circle (1cm and 2.25cm);
\draw[Red!50!Sepia,fill=Gold!50] (0,1.15) circle (.8cm and 1.75cm);
\draw[Sepia,fill=Gold!80] (0,.75) circle (.6cm and 1.25cm);
\foreach \n in {0,1,...,30} {
 \pgfmathsetmacro\s{2*exp(-\n/7)}
 \pgfmathsetmacro\c{int(50*\s)}
 \pgfmathsetmacro\cc{int(100-2*\n)}
 \draw (0,{3.95*(1-exp(-\n/5))}) node[Black!\c!Sepia!\cc,dot,scale=\s] (p\n) {};
 \ifthenelse{\n<11}{\draw (p\n) node[right=\s*1pt,scale=.5*\s] {$\n$};}{}
    }
\end{tikzpicture}
\end{wrapfigure}
Let me illustrate this approach using a very simple potentialist system that I refer to as \emph{Aristotelian} arithmetic potentialism, which might also be called \emph{initial-segment} potentialism, since the possible worlds consist of the truncations $\N\restrict n$ of the standard model of arithmetic. Each such world exists in the potentialist system of all such worlds, able to be extended to larger worlds as needed or desired, but no world has all the numbers.

\enlargethispage{30pt}%
What are the modal validities of this vision of potentialism? Well, worlds do not satisfy the assertion ``\emph{every number has a successor},'' written in symbols
 $$\forall a\,\exists b\ b=a+1$$
since the current largest number does not yet have a successor. But the worlds do satisfy the assertion ``\emph{necessarily, every number possibly has a successor},'' which can be expressed in the modal language as
 $$\necessary\forall a\, \possible\exists b\ b=a+1.$$
They also satisfy ``\emph{necessarily, every two numbers possibly have a sum and a product}''
 $$\necessary\forall a,b\ \possible\exists c,d\ (c=a+b\wedge d=a\cdot b),$$
and ``\emph{necessarily, for every number, possibly, there is a prime above it},'' which can also be expressed in the modal language.\goodbreak


Aristotelian potentialism is an instance of what I call \emph{linear inevitabilism}, in that the possible worlds are linearly ordered. In linear inevitabilism, the manner in which possibility unfolds through the worlds is always the same. The worlds are building up always in the same convergent manner to a unique limit structure, even if that limit structure is not a world in the potentialist system.

Potentialist systems with linear inevitabilism will always validate certain modal principles, such as the fact that every possibly necessary assertion is also necessarily possible.
 $$\possible\necessary\varphi\to\necessary\possible\varphi$$
This is the axiom (.2), the principal axiom of the modal logic S4.2, which is sound and complete for all convergent Kripke frames. Linear systems also validate a stronger principle, which asserts that if two assertions are both possible, then it is possible that one of them is true, while the other is possible.
$$\possible\varphi\wedge\possible\psi\to\left[ \possible(\varphi\wedge\possible\psi)\vee\possible(\psi\wedge\possible\varphi)\right]$$
This is the axiom (.3), the principal axiom of S4.3, which is sound and complete for all linearly ordered Kripke frames. So S4.3 is valid for Aristotelian arithmetic potentialism.


But let me next exhibit an alternative kind of arithmetic potentialism, which is not linear. For this system, we shall not take the numbers in their usual numerical order. To see why this might make some sense, consider the number googolplex, which is $10^{10^{100}}$. How does the ontological status of this number compare with the smaller numbers? Well, a typical smaller number has a googol digits. If those digits were essentially random, then to specify the number we would basically have to specify each and every digit separately---there wouldn't be a more efficient way to do it. But if one were able to recite, say, a million digits every second, then even if you did this from the beginning of time at the big bang, you will have recited only a tiny fraction of the digits of this typical number below a googolplex; there just isn't enough time to recite a googol digits in that duration. In particular, there is no reasonable sense in which a human being can hold a typical such number below a googolplex as an object of thought in mind. But meanwhile, we do seem able in contrast to hold a googolplex as an object of thought, for we can say quite a bit about it, such as the fact that it is even and not a multiple of 13 and so forth. A potentialist might reasonably want to say, therefore, that a googolplex might come into existence some time before most of the smaller numbers, which cannot be reasonably held as objects of thought. In short, this is a picture of arithmetic potentialism in which the realms of feasibility are not always downward closed.


\begin{wrapfigure}[11]{r}{.2\textwidth}\vskip-4ex\hfill
\begin{tikzpicture}[scale=1]
  \draw[use Hobby shortcut,closed=true,Sepia,fill=Yellow!20] (0,-1)..(-1,1)..(0,3)..(1,1);
  \draw[use Hobby shortcut,closed=true,red,fill=Tomato,fill opacity=.2] (0,-.7)..(-.7,1)..(0,2)..(.3,1.5)..(-.2,1)..(.6,.4);
  \draw[use Hobby shortcut,closed=true,Blue,fill=Orchid,fill opacity=.2] (0,-.8)..(-.6,.3)..(.3,1)..(0,2.5)..(.7,1);
  \draw[Sepia,thick,fill=Gold!80] (0,0) circle (.5cm) node {$u$};
  \draw[Red] (-.45,.85) node[scale=.8] {$v_0$};
  \draw[Blue] (.5,1.25) node[scale=.8] {$v_1$};
  \draw[Sepia] (0,2.8) node[scale=.8] {$w$};
  \draw (0,3) node[above] {$\vdots$};
\end{tikzpicture}
\end{wrapfigure}
\enlargethispage{20pt}%
So let me propose to consider \emph{arbitrary set} arithmetic potentialism, where the possible worlds consist of any finite set of numbers. In this potentialist conception, sometimes very large numbers come into actuality before some smaller numbers, or otherwise, and so it is genuinely nonlinear. This is a quite reasonable take on arithmetic potentialism, but it is fundamentally different from the Aristotelian conception and it will exhibit different modal validities. For example, modal axiom (.3) will not be valid for arbitrary set potentialism, although axiom (.2) will be valid, since the system is convergent in the sense that any two worlds are subsumed by some larger common extension.\goodbreak


This is therefore a particular instance of convergent potentialism, where the worlds are not necessarily linear ordered, but we have amalgamation---any two worlds we can access are subsumed by some yet further world. In this case, we still get the validity of axiom (.2) $\possible\necessary\varphi\to\necessary\possible\varphi$, but not always (.3) $\possible\varphi\wedge\possible\psi\to\left[ \possible(\varphi\wedge\possible\psi)\vee\possible(\psi\wedge\possible\varphi)\right]$, if the system isn't linear.

Every convergent potentialist system $\mathcal{W}$ converges to a limit model $W$, and indeed we need only that every world in $\mathcal{W}$ is a substructure of $W$ and every world in $\mathcal{W}$ can be extended so as to accommodate any additional desired individual of $W$. Both Aristotelian and arbitrary-set arithmetic potentialism are convergent to the standard model of arithmetic $\N$. More generally, examples of convergent potentialist systems arise if one takes the finite (or finitely generated) substructures of a given structure or the countable substructures of a fixed uncountable structure. For example, taking the countable transitive sets as a potentialist system of set theory is convergent to the set-theoretic universe $\HC$, consisting of the hereditarily countable sets, which is a model of $\ZFCm$ plus the countabilism axiom.

For any convergent potentialist system $\mathcal{W}$, we may undertake the \emph{potentialist translation}, translating actualist assertions $\psi$ to their modal counterparts $\psi^{\possible}$ simply by replacing all instances of $\exists x$ with $\possible\exists x$ and all instances of $\forall x$ with $\necessary\forall x$. The main point, appearing as theorem 1 in \cite{HamkinsLinnebo2022:Modal-logic-of-set-theoretic-potentialism}, is that if a potentialist system $\mathcal W$ of possible worlds converges to limit model $W$, then for any possible world $M\in\mathcal{W}$, we have
 $$W\satisfies \psi\quad\text{ if and only if }\quad M\satisfies_{\mathcal W} \psi^{\possible}.$$
In this way, actualist truth in the limit model of the system reduces to potentialist truth at the particular worlds of the potentialist system. The potentialist translation thus reveals a fundamental aspect of convergent potentialism. Namely, whilst using only the potentialist ontology, the convergent potentialist can thus give a full account of the actualist objects, the actualist structure, and even actualist truth. The actualist objects, after all, are exactly the possible objects of the potentialist system, and since this is a convergent system, any atomic structure they might come to have is inevitable and necessary---all the possible worlds will agree on the atomic structure of the objects. Finally, truth in the limit model is revealed in the modal system by the potentialist translation. In short, the potentialist denies the limit model exists, yet seems nevertheless to know everything that one would want to know about it.

%
To my way of thinking, this situation showcases not the power of the potentialist perspective, but rather how little there often is at stake in the dispute between potentialism and actualism, when the potentialist vision on offer is a form of convergent potentialism. In short, I say, convergent potentialism is \emph{implicitly actualist}. It is implicitly actualist because the convergent potentialist can give a completely clear account of the actualist model, interpreting the actualist model entirely within the potentialist ontology---there is nothing missing. For the convergent potentialist to deny the actualist model is something like accepting the real number field $\R$ but rejecting the complex field $\C$, even though $\C$ is interpretable in $\R$ using pairs of real numbers.



\newpage
\begin{wrapfigure}[12]{r}{.25\textwidth}\vskip-2ex\hfill
\begin{tikzpicture}[scale=.45,xscale=.7,line join=bevel]
\node (0) at (0,0) {};
\node (a) at (-1,4) {};
\node (b) at (1,4) {};
\node (c) at (-4,7.1) {};
\node (d) at (-1,7.1) {};
\node (e) at (0,6.2) {};
\node (f) at (2.5,6.2) {};
\node (g) at (2.1,9.5) {};
\node (h) at (5.1,9.5) {};
\node (i) at (-2,10) {};
\node (j) at (1.5,10) {};
 \draw[fill=blue,fill opacity=.16,thin] (0.center) -- (a.center) -- node[below,opacity=1,scale=.7] {$M$} (b.center) -- cycle;
 \draw[fill=blue,opacity=.9,fill opacity=.08] (0.center) -- (a.center) to[out=104,in=-60] (c.center) -- node[below,opacity=.9,scale=.6] {$M_0$} (d.center) to[out=-85,in=76] (b.center) -- cycle;
 \draw[fill=red!80!yellow,opacity=.8,fill opacity=.24,thin] (0.center) -- (a.center) to[out=104,in=-135] (e.center) -- node[below,opacity=.9,scale=.6] {$M_1$} (f.center) to[out=-135,in=76] (b.center) -- cycle;
 \draw[fill=red,opacity=.7,dash pattern={on 8pt off 1pt},fill opacity=.08] (e.center) to[out=45,in=-84] (g.center) -- node[below,opacity=.8,scale=.6] {$M_{11}$} (h.center) to[out=-125,in=45] (f.center);
 \draw[fill=yellow,opacity=.5,dotted,fill opacity=.08] (e.center) to[out=45,in=-75] (i.center) -- node[below,opacity=.8,scale=.6] {$M_{10}$} (j.center) to[out=-104,in=45] (f.center);
\end{tikzpicture}
\end{wrapfigure}
Meanwhile, there are other more radical forms of potentialism that are not convergent and that do not support the potentialist translation. There are forms of potentialism for which what is possible or necessary may depend on what's already happened. As some objects become actual, they may close off what used to be alternative possibilities. If a computation is revealed to converge with output $0$, for example, it will never subsequently have output $1$, even if that had been possible and indeed possibly necessary from an earlier world. Such is the nature of several of the forms of potentialism I have investigated in arithmetic and set theory.


Exactly this situation, for example, is realized in the potentialist system consisting of the models of \PA\ under end-extension. In this system, we view each (often nonstandard) model of \PA\ as a realm of feasibility, and moving to a larger realm is to invoke the possibility operator. What I proved is that if we take the nonstandard models of \PA\ as realms of feasibility, we achieve the radical-branching form of potentialism.

\begin{theorem}[Hamkins \cite{Hamkins:The-modal-logic-of-arithmetic-potentialism}]
Arithmetic end-extensional potentialism with the models of \PA\ has radical branching, and the modal validities are exactly only S4.
\end{theorem}

The proof makes use the universal algorithm, which is a Turing machine program that always enumerates a finite list of numbers, but in any model of \PA, you can make the next number on the sequence be whatever you would like by extending the model to a suitable end-extension. The nature of the numbers on the sequence produced by the universal algorithm, therefore, have a radical branching behavior, causing the entire potentialist system to exhibit radical branching for its modal possibilities.

The same phenomenon arises in set-theoretic potentialism, revealed in a sequence of papers. Woodin and I proved \cite{HamkinsWoodin:The-universal-finite-set} proved that the countable models of \ZFC\ under top-extensional potentialism exhibit radical branching, validating exactly S4. My student Julia (Kameryn) Williams and I proved \cite{HamkinsWilliams2021:The-universal-finite-sequence} that the countable models of \ZFC\ under end-extensional potentialism exhibit radical branching, validating exactly S4. In \cite{Hamkins2024:Every-countable-model-of-arithmetic-or-set-theory-has-a-pointwise-definable-end-extension} I used a further generalization to provide a flexible method of showing that every countable model of arithmetic or set theory has a pointwise definable end extension.


In this essay I have described several different varieties of potentialism, depending on the degree to which there is an inevitable coherency amongst the possible worlds. The tightest form occurs in a linear potentialist system, where possibility has the character of linear inevitability and generally validates the modal theory S4.3. Slightly more flexible is the nonlinear but convergent forms of potentialism, in which the worlds of the potentialist system converge to an actualist limit model. Because the objects, structure, and truths of this limit model are revealed in the potentialist system via the potentialist translation, I regard such potentialist conceptions as implicitly actualist. To my way of thinking, a truly potentialist ontology must allow for the radical branching case of mathematical possibility, for which we may be in a particular realm, but where the nature of possibility and necessity is itself contingent.

$$
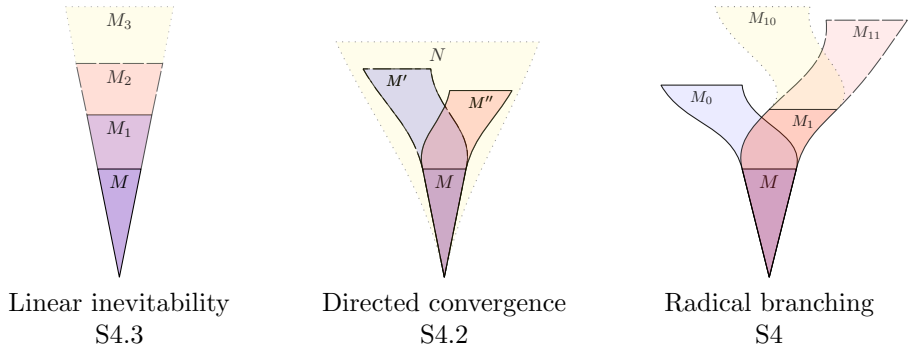
\begin{tikzpicture}[scale=.4,yscale=.9,xscale=.9,line join=bevel]
\node (0) at (0,0) {};
\node (a) at (0,4) {};
\node (b) at (0,6) {};
\node (c) at (0,7.9) {};
\node (d) at (0,10) {};
 \draw[fill=blue!67!red,fill opacity=.32,thin] (0.center) -- ([xshift=.8cm]a.center) -- node[below,opacity=1,scale=.7] {$M$} ([xshift=-.8cm]a.center) -- cycle;
 \draw[fill=blue!50!red,opacity=.7,fill opacity=.24,thin] ([xshift=.8cm]a.center) -- ([xshift=1.2cm]b.center) -- node[below,opacity=.9,scale=.7] {$M_1$} ([xshift=-1.2cm]b.center) -- ([xshift=-.8cm]a.center);
 \draw[fill=red!80!yellow,dash pattern={on 8pt off 1pt},opacity=.6,fill opacity=.16,thin] ([xshift=1.2cm]b.center) -- ([xshift=1.6cm]c.center) -- node[below,opacity=.8,scale=.7] {$M_2$} ([xshift=-1.6cm]c.center) -- ([xshift=-1.2cm]b.center);
 \draw[fill=yellow,dotted,opacity=.3,fill opacity=.08] ([xshift=1.6cm]c.center) -- ([xshift=2cm]d.center) -- node[below,opacity=.8,scale=.7] {$M_3$} ([xshift=-2cm]d.center) -- ([xshift=-1.6cm]c.center);
 \node[align=center] at (0,-1.5) {Linear inevitability\\ \theoryf{S4.3}};

\begin{scope}[shift={(12,0)}]
\node (0) at (0,0) {};
\node (a) at (-.8,4) {};
\node (b) at (.8,4) {};
\node (c) at (-3,7.7) {};
\node (d) at (-.5,7.7) {};
\node (e) at (.2,6.9) {};
\node (f) at (2.5,6.9) {};
 \draw[fill=blue!67!red,fill opacity=.24,thin] (0.center) -- (a.center) -- node[below,opacity=1,scale=.7] {$M$} (b.center) -- cycle;
 \draw[fill=blue,opacity=.8,dash pattern={on 7pt off 1pt}, fill opacity=.08] (a.center) to[out=104,in=-60] (c.center) -- node[below,opacity=.9,scale=.6] {$M'$} (d.center) to[out=-85,in=76] (b.center);
 \draw[fill=blue,opacity=.8,dash pattern={on 8pt off 2pt}, fill opacity=.08] (0.center) -- (a.center) to[out=104,in=-60] (c.center) -- node[below,opacity=.9,scale=.6] {$M'$} (d.center) to[out=-85,in=76] (b.center) -- cycle;
 \draw[fill=red,opacity=.7,fill opacity=.08,thin] (a.center) to[out=104,in=-100] (e.center) -- node[below,opacity=.8,scale=.6] {$M''$} (f.center) to[out=-120,in=76] (b.center);
 \draw[fill=red,opacity=.6,fill opacity=.08,thin] (0.center) -- (a.center) to[out=104,in=-100] (e.center) -- node[below,opacity=.8,scale=.6] {$M''$} (f.center) to[out=-120,in=76] (b.center) -- cycle;
 \draw[fill=yellow,dotted,opacity=.4,fill opacity=.08] (0.center) to[out=102,in=-65] ([shift={(-1,1)}]c.center) -- node[below,opacity=.8,scale=.7] {$N$} ([shift={(4,1)}]d.center) to[out=-110,in=78] cycle;
 \node[align=center] at (0,-1.5)  {Directed convergence\\ \theoryf{S4.2}};
\end{scope}

\begin{scope}[shift={(24,0)}]
\node (0) at (0,0) {};
\node (a) at (-1,4) {};
\node (b) at (1,4) {};
\node (c) at (-4,7.1) {};
\node (d) at (-1,7.1) {};
\node (e) at (0,6.2) {};
\node (f) at (2.5,6.2) {};
\node (g) at (2.1,9.5) {};
\node (h) at (5.1,9.5) {};
\node (i) at (-2,10) {};
\node (j) at (1.5,10) {};
 \draw[fill=blue,fill opacity=.16,thin] (0.center) -- (a.center) -- node[below,opacity=1,scale=.7] {$M$} (b.center) -- cycle;
 \draw[fill=blue,opacity=.9,fill opacity=.08] (0.center) -- (a.center) to[out=104,in=-60] (c.center) -- node[below,opacity=.9,scale=.6] {$M_0$} (d.center) to[out=-85,in=76] (b.center) -- cycle;
 \draw[fill=red!80!yellow,opacity=.8,fill opacity=.24,thin] (0.center) -- (a.center) to[out=104,in=-135] (e.center) -- node[below,opacity=.9,scale=.6] {$M_1$} (f.center) to[out=-135,in=76] (b.center) -- cycle;
 \draw[fill=red,opacity=.7,dash pattern={on 8pt off 1pt},fill opacity=.08] (e.center) to[out=45,in=-84] (g.center) -- node[below,opacity=.8,scale=.6] {$M_{11}$} (h.center) to[out=-125,in=45] (f.center);
 \draw[fill=yellow,opacity=.5,dotted,fill opacity=.08] (e.center) to[out=45,in=-75] (i.center) -- node[below,opacity=.8,scale=.6] {$M_{10}$} (j.center) to[out=-104,in=45] (f.center);
 \node[align=center] at (0,-1.5)  {Radical branching\\ \theoryf{S4}};
\end{scope}
\end{tikzpicture}$$

In conclusion, I should like to propose a vision of ultrafinitism engaging more directly with potentialist conceptions. Ultrafinitism is naturally concerned with the concept of feasibility, but this is a concept that clearly admits different natures and instantiations at different scales, and in this ultrafinitism will have a potentialist modal character.

We have already seen how the models of finite arithmetic \FA\ exhibit this potentialist nature in their interpretations to taller structures. And it is surely in keeping with the ultrafinitist attitude that no one realm of feasibility would have the whole story about how all the rest of the possible realms will unfold, so we should not think that we are limited only to truncated models of that minimal extension or to the resulting model $M^*$ of bounded induction I$\Delta_0$. Rather, the ultrafinitist discussion should be about: what is the nature of a realm of feasibility? If we write down an arithmetic theory that we expect to be true in any such realm, then there is a natural potentialist system available---just consider all the models of that theory, each accessing all its extensions to further models of that theory.

My theorem above shows that in the case of \PA, we get an elegant potentialist theory, validating exactly the modal theory S4. A similar analysis works for I$\Sigma_1$, which supports an analogue of the universal algorithm, and perhaps this is also true for I$\Delta_0+\text{exp}$, in which exponentiation is total. But the proof breaks down in weaker theories such as I$\Delta_0$ itself or \FA.

So there is a natural proposal here for research in modal potentialist ultrafinitism.

\begin{questions*}\
 \begin{enumerate}
  \item What are the modal validities of the potentialist system of all models of \FA\ under end-extension?
  \item What are the modal validities of the potentialist system of all models of I$\Delta_0$ under end-extension?
 \end{enumerate}
\end{questions*}

\noindent
(For the models of \FA\ under end-extension, one should re-interpret the constant $N$ for the largest number in each world separately, in effect taking this constant de dicto rather than de re; or alternatively, since the largest number of any world is definable in that world, simply omit this constant from the language.) Perhaps the theory $\necessary\FA+\text{S4}$ expresses a core ultrafinitist potentialist commitment, asserting that \FA\ holds in all worlds of a potentialist context. What is the underlying modal logic?

\printbibliography

\end{document}